\documentclass[10pt]{article}
\usepackage{geometry}                
\usepackage{graphicx}
\usepackage{subfig}
\usepackage{amssymb}
\usepackage{amsmath}
\usepackage{hyperref}
\usepackage{epstopdf}
\usepackage{caption}
\newcommand{\vecb}{\boldsymbol{\mathbf{b}}}

\newcommand{\vece}{\boldsymbol{\mathbf{e}}}

\newcommand{\veck}{\boldsymbol{\mathbf{k}}}

\newcommand{\vecm}{\boldsymbol{\mathbf{m}}}

\newcommand{\vecw}{\boldsymbol{\mathbf{w}}}
\newcommand{\vecx}{\boldsymbol{\mathbf{x}}}
\newcommand{\vecy}{\boldsymbol{\mathbf{y}}}

\newcommand{\matA}{\boldsymbol{\mathbf{A}}}
\newcommand{\matB}{\boldsymbol{\mathbf{B}}}

\newcommand{\matH}{\boldsymbol{\mathbf{H}}}
\newcommand{\matI}{\boldsymbol{\mathbf{I}}}

\newcommand{\matK}{\boldsymbol{\mathbf{K}}}

\newcommand{\matS}{\boldsymbol{\mathbf{S}}}

\newcommand{\matU}{\boldsymbol{\mathbf{U}}}

\newcommand{\matW}{\boldsymbol{\mathbf{W}}}
\newcommand{\matX}{\boldsymbol{\mathbf{X}}}
\newcommand{\matY}{\boldsymbol{\mathbf{Y}}}

\newcommand{\matzero}{\boldsymbol{\mathbf{0}}}

\newcommand{\tv}{\boldsymbol{\theta}}

\newcommand{\xv}{\boldsymbol{\xi}}

\newcommand{\bv}{\boldsymbol{\beta}}

\title{Covariate dimension reduction for survival data\\via the Gaussian process latent variable model}
\author{James E. Barrett\footnote{Contact: regmjeb@kcl.ac.uk}\,\,\, and Anthony C. C. Coolen\\
\small{Institute of Mathematical and Molecular Biomedicine, King's College London}}
\date{November 2015}

\begin{document}
\maketitle
\begin{abstract}
The analysis of high dimensional survival data is challenging, primarily due to the problem of overfitting which occurs when spurious relationships are inferred from data that subsequently fail to exist in test data. Here we propose a novel method of extracting a low dimensional representation of covariates in survival data by combining the popular Gaussian Process Latent Variable Model (GPLVM) with a Weibull Proportional Hazards Model (WPHM). The combined model offers a flexible non-linear probabilistic method of detecting and extracting any intrinsic low dimensional structure from high dimensional data. By reducing the covariate dimension we aim to diminish the risk of overfitting and increase the robustness and accuracy with which we infer relationships between covariates and survival outcomes. In addition, we can simultaneously combine information from multiple data sources by expressing multiple datasets in terms of the same low dimensional space. We present results from several simulation studies that illustrate a reduction in overfitting and an increase in predictive performance, as well as successful detection of intrinsic dimensionality. We provide evidence that it is advantageous to combine dimensionality reduction with survival outcomes rather than performing unsupervised dimensionality reduction on its own. Finally, we use our model to analyse experimental gene expression data and detect and extract a low dimensional representation that allows us to distinguish high and low risk groups with superior accuracy compared to doing regression on the original high dimensional data. 
\end{abstract}
\section{Introduction}
\label{sec:intro}

High dimensional data are increasingly common in biomedical research. For instance, current experimental techniques can acquire tens of thousands of gene expression measurements or hundreds of thousands of SNP (single nucleotide polymorphism) data. Automated image analysis software can generate hundred or thousands of parameters from biomedical images obtained using various imaging platforms. The analysis of high dimensional data is a challenging problem and this is also the case with high dimensional survival data, where in addition to covariates we also measure the time until an event of interest. One of the main difficulties is \emph{overfitting} which occurs when a model fits training data very well but fails to generalise to test data. This happens when a model fits to noise and struggles to detect genuine structure in the data. The greater the dimension of the data compared to the number of samples the more difficult it becomes to extract meaningful relationships between the covariates and outcomes. Applying traditional methods such as a Cox Proportional Hazards Model \cite{COX72} is problematic as the regression coefficients are not uniquely determined when the number of covariates ($d$) exceeds the number of samples ($N$) \cite{WIT10}.

Strategies for tackling high dimensional data can be divided into two broad classes. \emph{Supervised} methods take into account the survival outcomes. For example, feature selection aims to select a subset of the covariates that are relevant either by doing a univariate analysis on each covariate and selecting the most significant \cite{WIT10} or performing $L_1$ or $L_2$-penalised regression with a Cox model \cite{GOE10, PARK07, SOHN09}. Random forests \cite{ISH11} and elastic nets \cite{ENG09} have also been proposed for feature selection with survival data. These approaches are suitable when the goal is to establish associations between covariates and survival outcomes.

Alternative \emph{unsupervised} dimensionality reduction methods attempt to represent the information in a high dimensional dataset in a lower dimensional space. The idea is that there will in general be some redundancy between high dimensional covariates, and that by eliminating this redundancy we can achieve a more parsimonious representation of the data. Since the ratio of covariates to samples is now smaller, the risk of detecting spurious relationships is reduced. These approaches are appropriate when the goal is to make predictions for new individuals since overfitting will hopefully be diminished and consequently predictive accuracy will increase. A drawback is that the impact a particular covariate has on the survival outcomes may not be easy to interpret (since the low dimensional representation may be a complicated combination of the high dimensional covariates). For an excellent overview of survival analysis with high dimensional data see \cite{WIT10}.

One approach to dimensionality reduction is via latent variable models which attempt to represent the information contained in a high dimensional dataset in terms of a smaller number of latent variables. In this paper we extend the popular Gaussian Process Latent Variable Model \cite{LAW05} (GPLVM) to incorporate survival outcomes. The GPLVM is a flexible probabilistic non-linear dimensionality reduction method. The model assumes that the high dimensional covariates can be written as a stochastic function of the latent variables and assumes a Gaussian process (GP) prior over that function. By choosing different kernel functions in the GP prior various types of non-linear mappings can be specified between the low and high dimensional spaces. The latent variables are unknown and must be inferred from the data.

The simplest case consists of a linear mapping from the low to high dimensional spaces (corresponding to a linear GP kernel). It was shown in \cite{LAW05} that the maximum a posteriori (MAP) solution for the latent variables is equivalent to performing Principal Component Analysis (PCA) and retaining the first $q$ principal components (where $q$ is the number of latent variables). We can intuitively regard the GPLVM as a non-linear probabilistic generalisation of PCA.

A drawback of the original GPLVM is its computational complexity. This prompted the subsequent application of sparse GP regression methods to the GPLVM \cite{LAW07}. Recent advances in variational sparse Gaussian Process (GP) regression \cite{TIT09} have also been successfully applied \cite{LAW10}. A variational lower bound on the marginal likelihood was constructed which can then be optimised with respect to the variational parameters and model hyperparameters. A detailed overview can be found in \cite{GAl14b}. It is also possible to infer what the intrinsic dimensionality of the latent variable space is using this method (that is, how many latent variables are required to explain the observed data).

Another use of the GPLVM has been to combine multiple sources of data by simultaneously expressing several datasets in terms of the same latent variables. The idea is that overlapping structure can be easily captured by shared latent variables. This has been developed in the shared-GPLVM \cite{SHO06,EK07}. That work was further extended to allow each dataset to have a separate set of latent variables that would account for information unique to each source \cite{EK08,DAM12}. There have also been extensions of the model to include `output' information. In the discriminative-GPLVM \cite{URT07,ELEF13} class labels are incorporated and a low dimensional embedding is extracted that attempts to minimise within-class variance and maximise between-class variance. The supervised-GPLVM \cite{GAO11} includes continuous output variables.

The main advance in this paper is to incorporate (possibly censored) survival outcomes into the GPLVM by combing the GPLVM with a Weibull proportional hazards model (WPHM). The latent variables now attempt to simultaneously capture structure contained in both the high dimensional covariates and the survival outcomes. By combining both sources of information we hope to infer a low dimensional representation that captures not just the low dimensional structure of the covariates but also the relationship between covariates and outcomes. By connecting the covariates to survival outcomes via the low dimensional latent variable space we are limiting the degrees of freedom the model has, and thereby reduce the risk of overfitting.

Recently the GPLVM has been applied to facial expression recognition \cite{ELEF13} which provides a useful analogy for the model proposed here. Images of a subject's face were taken from two different angles and these images were regarded as two different datasets. Both datasets are expressed in terms of the same latent variables since both images are of the same facial expression, but from different angles. In our model we can think of the latent variable space as representing the underlying biological processes we are interested in. Each observed dataset provides a different `view' or `perspective' onto those processes. For example, we may acquire gene expression data from cancer patients. Additionally, we might generate relevant parameters from imaging their tumours. If the gene expression data and imaging parameters are attempting to characterise the cancer in different ways then it is reasonable to propose that they offer two different but complementary `views' of the underlying tumour. In addition, if the survival outcomes are driven by those underlying biology then they provide yet another `view' and thus it is desirable to represent all of the observed data in terms of a shared low dimensional structure.

We overcome some technical issues to construct a Laplace approximation of the marginal likelihood. The Laplace approximation is straightforward to apply to complicated likelihoods involving WPHM terms and is used for the purposes of hyperparameter optimisation and model selection. We compare the model likelihood corresponding to different choices of $q$ (number of latent variables) in order to determine what the optimal dimensional of the latent variable space is. In our model we also allow for multiple datasets via a set of shared latent variables as in \cite{SHO06,EK08}.

We conduct several numerical simulations to study the effects of overfitting due to high dimensionality and examine the performance of the combined GPLVM-WPHM to detect and extract low dimensional structure under various conditions and finally apply the new model to gene expression data from the breast cancer METABRIC dataset. One of the goals in the METABRIC study was to identify potential new gene signatures that were associated with clinical outcome (overall survival or progression free survival). The risk of overfitting is high given the large number of genes and it is therefore desirable to reduce the dimension of the dataset while searching for potential associations between genes and clinical outcome. We show that the GPLVM-WPHM can achieve a predictive accuracy (measured using $k$-fold cross validation) that is considerably higher than using the WPHM alone, thereby achieving an advantage of practical benefit.

The rest of this paper is structured as follows. In Section \ref{sec:def} we define define the GPLVM and the WPHM separately before defining the combined GPLVM-WPHM. We provide details of the Laplace approximations, inference of parameters and hyperparameters, and how to make survival predictions for new individuals. In Section \ref{sec:sim} we present results from simulation studies and real data, and we finish with a discussion in Section \ref{sec:disc}.

\section{Model definition}
\label{sec:def}

\subsection{The Gaussian process latent variable model (GPLVM)}

\begin{figure}
\centering
\includegraphics[scale=0.8]{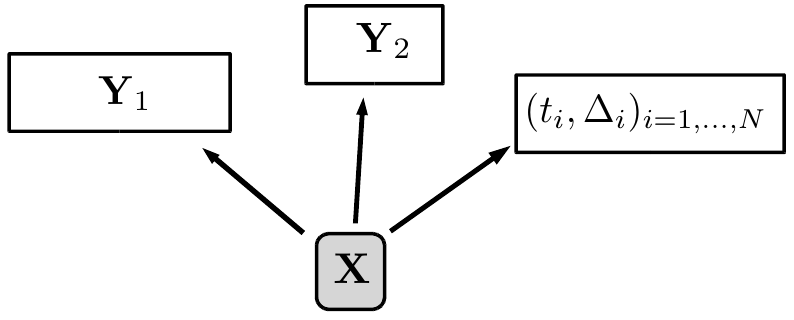}
\caption{Schematic diagram of the combined GPLVM-WPHM. Here two high dimensional datasets $\matY_1\in\mathbb{R}^{N\times d_1}$ and $\matY_2\in\mathbb{R}^{N\times d_2}$ and the survival data $(t_i,\Delta_i)$ for $i=1,\ldots,N$ are represented by the same underlying set of $q$ latent variables $\matX\in\mathbb{R}^{N\times q}$. The total number of patients is $N$, $d_1$ and $d_2$ are the dimensions of the observed datasets respectively, $t_i$ is the time to event for individual $i$ and $\Delta_i\in(0,1)$ is an indicator variables that tells us whether that individual was censored or experienced the primary event. The latent variable representation will attempt to represent information from all the sources of data in a more parsimonious form. By `squashing' the data into a smaller number of latent variables we aim reduce the risk of overfitting.}
\end{figure}

We consider $S$ observed datasets $\matY_1,\ldots,\matY_S$, each with $N$ rows which correspond to individuals, and $d_1,\ldots,d_S$ columns, respectively, which correspond to the covariates. For example, $\matY_1$ could be gene expression data, $\matY_2$ imaging parameters and so forth. It is assumed that the rows in each dataset correspond to the same individuals. We assume each individual $i$ can be represented by a low dimensional vector of latent variables $\vecx_i\in\mathbb{R}^q$ via
\begin{equation}
y_{i\mu}^s = f^s_{\mu}(\vecx_i) + \xi_{i\mu}^s\quad\text{for $i=1,\ldots,N$}
\label{model_eqn}
\end{equation}
where $y_{i\mu}^s$ is the $\mu$th covariate value for individual $i$ in dataset $s$. The functions $f^s_{\mu}$ are unspecified and a Gaussian process (GP) prior is assumed. The same prior is used for each dimension $\mu$ of the covariate space $s$, but may differ with respect to $s$. The noise variables $\xi_{i\mu}^s$ are i.i.d Gaussian random variables with zero mean and variances $\beta^2_s$. It is assumed that $q<\min_s(d_s)$. It follows that the data likelihood is \cite{LAW05}
\begin{equation}
p(\{\matY\}|\matX,\{\beta^2\},\tv) =\prod_{s=1}^S\prod_{\mu=1}^{d_s} \frac{e^{-\frac{1}{2}\vecy^s_{:,\mu}\cdot\matK_s^{-1}\vecy^s_{:,\mu}}}{(2\pi)^{\frac{N}{2}}|\matK_s|^{\frac{1}{2}}}
\label{eq:gplvm:datalikelihood}
\end{equation}
where $\{\matY\}$ denotes the set of all observed datasets, $\matX\in\mathbb{R}^{N\times q}$ is the matrix of latent variables, $\{\beta\}$ and $\tv$ are hyperparameters which are defined below. The kernel matrix is $[\matK_s]_{ij} = k_s(\vecx_i,\vecx_j) + \beta_s^2\delta_{ij}$ where $k(.,.)$ is called the kernel function. The kernel function tells us how `similar' $\vecx_i$ and $\vecx_j$ are. If two individuals are `similar' in the latent space then they are also likely to be `similar' in the observed covariate space. The vector $\vecy^s_{:,\mu}$ denotes the $\mu$th column of $\matY_s$. Equation (\ref{eq:gplvm:datalikelihood}) is thus a product of $d_s$ Gaussian processes for each dataset which map the latent variables to each covariate in $\matY_s$. A Gaussian process prior \cite{RAS06} can be regarded as a prior over functions and is completely specified by its mean function, $m(\vecx_i) = \text{mean}(\vecx_i)$ (zero in this case), and the kernel function $k(\vecx_i,\vecx_j) = \text{cov}(\vecx_i,\vecx_j)$ where $\text{cov}(\vecx_i,\vecx_j)$ denotes the covariance between $\vecx_i$ and $\vecx_j$. The kernel functions considered in this paper are
\begin{equation}
\begin{array}{ll}
k(\vecx_i,\vecx_j) = \sigma\vecx_i\cdot\vecx_j & \text{linear,}\\
k(\vecx_i,\vecx_j) = \sigma(1+\vecx_i\cdot\vecx_j)^2& \text{polynomial (of second order),}\\
k(\vecx_i,\vecx_j) = \sigma \exp(-l(\vecx_i-\vecx_j)^2/2)& \text{squared exponential.}
\end{array}
\label{eq:kernels}
\end{equation}
In all three kernels the hyperparameter $\sigma$ controls the variance of the `outputs' (which in our case are the high dimensional covariates). The hyperparameter $l$ defines a characteristic length scale over which the values of the outputs vary.

Consider for the moment that $s=1$. With linear functions $f_{\mu}$ in (\ref{model_eqn}) we can write $\vecy_i=\matW\vecx_i+\xv_i$ where $\matW\in\mathbb{R}^{N\times q}$ is a matrix of mapping coefficients. As shown originally in \cite{LAW05} if we assume a Gaussian prior for $\matW$ and marginalise we obtain equation (\ref{eq:gplvm:datalikelihood}) with the linear kernel. The MAP solution for $\matX$ in this case is equivalent to performing principal component analysis and retaining the top $q$ principal components.
\subsection{Weibull proportional hazards model (WPHM)}

For each individual $i$, in addition to the covariates we observe an event time $t_i$ and an indicator variable where $\Delta_i=1$ means the primary event occurred first and $t_i$ is therefore the time until the primary event, whereas $\Delta_i=0$ indicates that that individual was right censored (censoring event times are assumed to be independent of the primary risk event times). In the WPHM the hazard rate for individual $i$ is
\begin{equation}
h_i(t|\vecx_i,\nu,\rho,\vecb) = \lambda_0(t)e^{\vecb\cdot\vecx_i},
\label{eq:wphm:hazard}
\end{equation}
where the base hazard rate is $\lambda_0(t) = (\nu/\rho)(t/\rho)^{\nu-1}$. The scale parameter $\rho$, shape parameter $\nu$ and regression coefficients $\vecb\in\mathbb{R}^q$ need to be inferred from the data. Note that in anticipation of combining the GPLVM with the WPHM the hazard rate is a function of the latent variables $\vecx_i$ and not the observed data. Denoting the survival data as $D=\{(t_1,\Delta_1),\ldots,(t_N,\delta_N)\}$ and the integrated base hazard rate as $\Lambda_0(t) = (t\rho)^{\nu}$, the data likelihood is
\begin{equation}
p(D|\matX,\vecb,\rho,\nu) = \prod_{i=1}^N [\lambda_0(t_i)e^{\vecb\cdot\vecx_i}]^{\Delta_i}\exp(-\Lambda_0(t_i)e^{\vecb\cdot\vecx_i}).
\end{equation}
Note that since this is a parametric model would be straightforward to include left, right or interval censored observations.
\subsection{The combined GPLVM-WPHM}
\label{modeldef}
We are now interested in a model where the high dimensional covariates and the survival outcomes are both related to the same latent variables. Using Bayes' theorem we can write the posterior density over the unknown parameters:
\begin{equation}
p(\matX,\vecb,\rho,\nu|\{\matY\},D,\tv,\{\beta\}) = \frac{p(\{\matY\},D|\matX,\vecb,\rho,\nu,\tv,\{\beta\})p(\matX)p(\vecb)p(\rho)p(\nu)}{p(\{\matY\},D|\{\beta\},\tv)}
\label{eq:joint_posterior}
\end{equation}
where
\begin{equation}
p(\{\matY\},D|\{\beta\},\tv) = \int \text{d}\matX\,\text{d}\vecb\,\text{d}\rho\,\text{d}\nu\, p(\{\matY\},D|\matX,\vecb,\rho,\nu,\tv,\{\beta\})p(\matX)p(\vecb)p(\rho)p(\nu).
\label{slvm:eq:marginal}
\end{equation}
We work with the negative log posterior $\mathcal{L}(\matX,\vecb,\rho,\nu;\{\beta\},\tv) = -N^{-1}\log p(\matX,\vecb,\rho,\nu|\{\matY\},D,\tv,\{\beta\})$ in practice. We now make a key assumption of conditional independence between the observed covariates and the survival data given the latent variables:
\begin{equation}
p(\{\matY\},D|\matX,\vecb,\rho,\nu,\tv,\{\beta\}) = p(\{\matY\}|\matX,\tv,\{\beta\})p(D|\matX,\vecb,\rho,\nu,\tv).
\end{equation}
The first term is given by the the GPLVM likelihood (\ref{eq:gplvm:datalikelihood}) and the second is given by the WPHM likelihood (\ref{eq:wphm:hazard}). Following the example of \cite[Section 2.2]{ib2001} we choose Gamma prior distributions for the scale and shape parameters:
\begin{equation}
p(\nu|\kappa_0,\alpha_0) = \frac{\nu^{\kappa_0-1}e^{-\nu/\alpha_0}}{\alpha_0^{\kappa_0}\Gamma(\kappa_0)}\quad\text{and}\quad p(\rho|\kappa_1,\alpha_1) = \frac{\rho^{\kappa_1-1}e^{-\rho/\alpha_1}}{\alpha_1^{\kappa_1}\Gamma(\kappa_1)}.
\label{eq:priors}
\end{equation}
For the regression parameters prior we chose $p(\vecb)=\mathcal{N}(\matzero,\sigma_0^{-2}\matI)$. We found that $(\kappa_0,\alpha_0) = (3,1),(\kappa_1,\alpha_1) = (3,6)$ and $\sigma_0 = 2$ to be satisfactory in practice (where we expect the event times are measured in years).

For the linear and polynomial kernel there is some redundancy between the overall scale controlled by $p(\matX)$ and the hyperparameter $\sigma$. We decided to fix $\sigma=1$ and use a flat improper prior $p(\matX) = 1$. The overall scale of $\matX$ is now naturally determined by the observed data (which will typically be normalised to unit variance and zero mean). In the case of the squared exponential kernel we imposed $\vecx_i\sim\mathcal{N}(\matzero,\sigma_1^{-2}\matI)$ with $\sigma_1=2$ for $i=1,\ldots,N$.
\subsection{Inference of parameters and hyperparamters}
\label{infer}

The latent variables, $\matX^*$ are determined by numerically solving $\matX^* = \min_{\matX}\mathcal{L}(\matX,\vecb,\rho,\nu)$ while holding $\vecb$, $\rho$ and $\nu$ fixed. This is followed by solving $(\vecb^*,\rho^*,\nu^*) = \min_{\vecb,\rho,\nu}\mathcal{L}(\matX^*,\vecb,\rho,\nu)$ where $\matX$ is fixed to its previously determined optimal value. This procedure is then repeated by alternately optimising with respect to one set of parameters while the others are fixed at their previously optimal values until a stable solution is converged upon. Further details of the implementation are given in Section \ref{sec:algorithm}. The posterior over hyperparameters is
\begin{equation}
p(\{\beta^2\},\tv|\{\matY\},D) = \frac{p(\{\matY\},D|\{\beta^2\},\tv)p(\{\beta^2\})p(\tv)}{\int\text{d}\{\beta^{\prime 2}\}\text{d}\tv'\,p(\{\matY\},D|\{\beta^{\prime 2}\},\tv')p(\{\beta^{\prime 2}\})p(\tv')}
\end{equation}
where the marginal likelihood $p(\{\matY_s\},D|\{\beta^2_s\},\tv)$ is defined by (\ref{slvm:eq:marginal}). Flat priors are assumed for the kernel parameters. The marginal likelihood involves an integral which is generally intractable both analytically and numerically. We therefore construct a Laplace approximation of the marginal likelihood. The negative hyperparameter log likelihood is defined as $-N^{-1}\log p(\{\beta^2\},\tv|\{\matY\})$ which in this case gives
\begin{equation}
\mathcal{L}_{hyp}(\{\beta^2\},\tv|\{\matY\})=\mathcal{L}(\matX^*,\vecb^*,\rho^*,\nu^*)-\frac{q}{2}\log2\pi + \frac{1}{2N}\log|N\matH(\{\beta^2\},\tv)|.
\label{eq:LLhyp}
\end{equation}
The Hessian matrix $\matH$ contains second order partial derivatives. Full details are given in the Supporting Information. Optimal hyperparameters are determined by numerically minimising (\ref{eq:LLhyp}). Note that each evaluation of $\mathcal{L}_{hyp}$ requires computing $(\matX^*,\vecb^*,\rho^*,\nu^*)$. This is computationally expensive although the search can be initialised to the optimal value from the previous evaluation of $\mathcal{L}_{hyp}$.

\subsection{Elimination of symmetries due to invariance under unitary transformations}
\label{sec:invariance}

A problem which arises with the Laplace approximation is that the second order partial derivatives are zero along certain directions in the $Nq$-dimensional posterior search space. These directions point along lines where the log likelihood is constant. This occurs due to the invariance of the log likelihood to rotation or reflection of the latent variables. A consequence of this is that the Hessian matrix is not guaranteed to be positive definite, and secondly the solution will not be unique (since any rotation or reflection is an equivalent solution). The problem can be eliminated by breaking the symmetry. We achieve this by constraining some of the degrees of freedom $\matX$ can take. This is fully described in the Supporting Information.

\subsection{Making predictions for new individuals}
\label{predictions}

When we observe a new individual with covariates $\vecy^*$ we wish to firstly infer their optimal position $\vecx^*$ in the latent variable space using the GPLVM and from there make a prediction of survival outcomes using the WPHM part of the model. To find the optimal location in the latent variable space we use the GP predictive distribution to write $p(\vecy^*|\vecx^*,\{\matY\},D) = \mathcal{N}(\vecm,\kappa^{-2}\matI)$ where $m_{\mu} = \veck\cdot\matK^{-1}\vecy_{\mu}$ and $\kappa^2 = k(\vecx^*,\vecx^*) - \veck\cdot\matK^{-1}\veck + \beta^2$. We maximise the posterior using gradient based methods (see Supporting Information for partial derivatives):
\begin{equation}
p(\vecx^*|\vecy_1^*,\ldots\vecy_S^*,\{\matY\},D) \propto p(\vecy_1^*|\vecx^*,\{\matY\},D)\cdots p(\vecy_S^*|\vecx^*,\{\matY\},D)p(\vecx^*).
\label{eq:pred}
\end{equation}
It is not necessary to make observations for each data source $s$ since the terms corresponding to missing data in (\ref{eq:pred}) can simply be omitted. Once we have determined $\vecx^*$ by maximising (\ref{eq:pred}) we use the WPHM part of the model to make survival predictions. We can generate a prediction of the event time, $t^*$, corresponding to $\vecx^*$ by numerically computing the mean of the corresponding event time density:
\begin{equation}
t^* = \left<t\right> = \int_0^{\infty}\text{d}s\,s\lambda_0(s)e^{\hat{\vecb}\cdot\vecx^*}\exp(-\Lambda_0(s)e^{\hat{\vecb}\cdot\vecx^*}).
\label{eq:etd}
\end{equation}
Note that the optimal values $\hat{\vecb}$ are used. Optimal values of $\nu$ and $\rho$ are used inside $\lambda_0(s)$ and $\Lambda_0(s)$. We can similarly compute the variance $\left<t^2\right>-\left<t\right>^2$ as a measure of the uncertainty associated with our prediction.

\subsection{Implementation}
\label{sec:algorithm}

Gradient based optimisation functions were used in a Matlab implementation. Partial derivatives are given in the Supporting Information. Initial values were set as $\vecb = \matzero$, $\rho=3$ and $\nu=10$. The initial values of $\matX$ are randomly generated from a Gaussian density with zero mean and covariance matrix equal to the identity matrix. In the case of the linear kernel function (\ref{eq:kernels}) it was shown in \cite{LAW05} that the GPLVM log likelihood has a single global minimum which corresponds to performing principal component analysis and retaining the top $q$ principal components. Experience suggests that in the GPLVM-WPHM the log likelihood still has a single minimum although this has not been proved. In the case of a non-linear kernel then there will exist multiple local minima. Several attempts are made to locate the global minimum, with each attempt starting from a different initial search point.

Software is available to download from the author's Mathwork's file exchange page. A model with $N=100$ and $q=2$ can be fitted in $<\mathcal{O}(10)$ mins on an intel i7 quad-core CPU. It is also possible to optimise the posterior (\ref{eq:joint_posterior}) with respect to $\{\beta\}$ and $\tv$ also. However, we found that the model will typically fail to infer the correct value of $q$ (larger values of $q$ have a higher posterior probability). The Laplace approximation was sufficient to penalise larger values of $q$ and the model performs well. As an alternative to the Laplace approximation one could compute the Bayes Information Criterion score for each value of $q$ and this may provide an acceptable level of performance at a reduced computational cost.

\section{Simulation studies}
\label{sec:sim}

\subsection{Generation of synthetic data}

\begin{figure}[t]
\centering
\includegraphics[scale = 0.8]{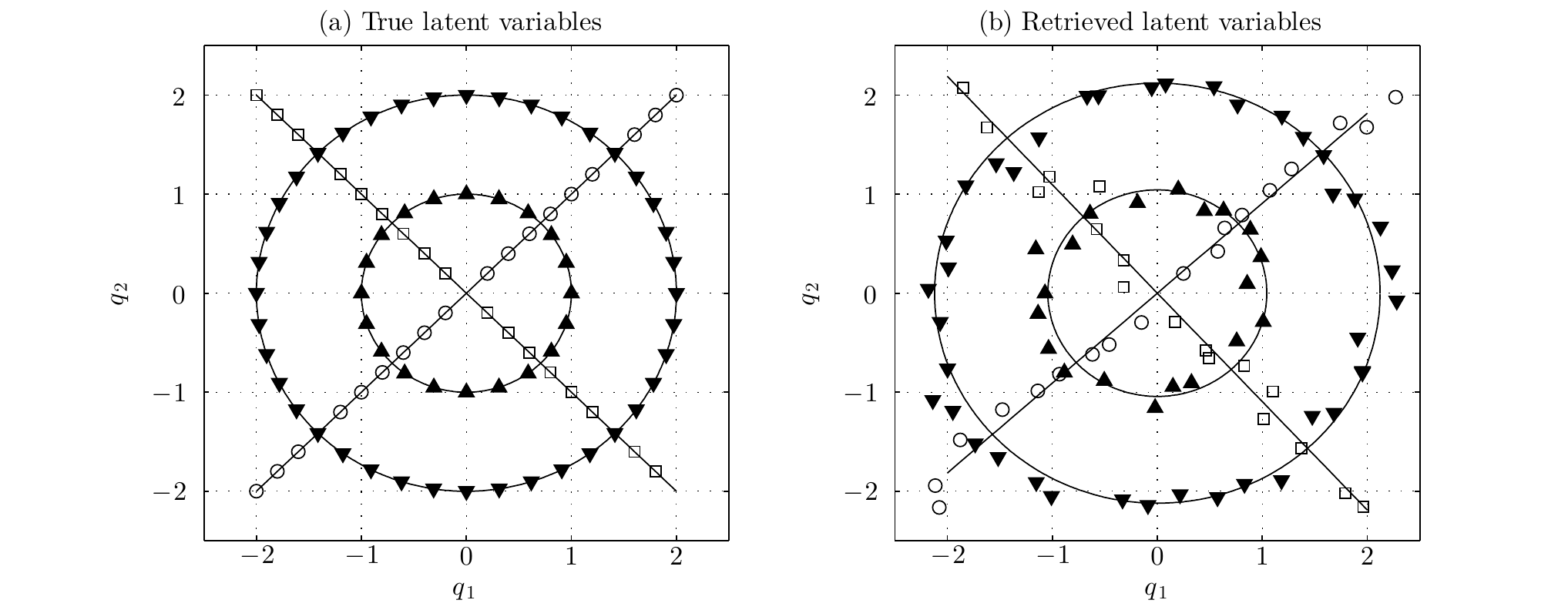}
\caption{On the left are the `true' synthetic latent variables that are projected into a high dimensional space with $d=10$ (using a linear kernel) according to (\ref{model_eqn}). The variance of the Gaussian noise was $\beta^2 = 0.1$. On the right are the latent variables retrieved by the GPLVM-WPHM model from those high dimensional data. This specific geometric pattern was used to allow for qualitative and quantitative comparison. Misalignment errors are $\mathcal{E}_{radial} = 0.0051$, $\mathcal{E}_{angular}= 0.0086$ and $\mathcal{E}_{linear}= 0.0288$.}
\label{special_X}
\end{figure}

To examine the behaviour of the model under different conditions we generate simulated data. We do this by first generating latent variables $\matX$ and from these generating high dimensional covariates $\matY$ and survival data $(t_i,\Delta_i)_{i=1,\ldots,N}$. The high dimensional data can be generated according to equation (\ref{eq:gplvm:datalikelihood}). Firstly the kernel matrix $\matK$ is computed (for certain chosen values of the hyperparameters $\tv$ and noise level $\beta$) and then for each dimension of $\matY$ we draw a random vector from a Gaussian process prior (which for fixed $N$ is simply a multivariate Gaussian density). This can be done for an arbitrary dimension $d$ and can be repeated to generate multiple datasets (with possibly different kernel functions and dimensions).

To generate survival outcomes we pick values for $\bv$, $\rho$, $\nu$ manually. Event times are generated using the inverse of the cumulative distribution which is  $C_i(t) = 1-\exp(-\Lambda_0(t)e^{\vecb\cdot\vecx_i})$ by generating random numbers $z\in[0,1]$ from a uniform density and obtaining the corresponding event time from the inverse cumulative distribution:
\begin{equation}
t_i = \rho\left(-e^{-\vecb\cdot\vecx_i}\log(1-z)\right)^{1/\nu}.
\end{equation}
Finally independent censoring is simulated by randomly selection a subset (about 10\%) of the individuals and generating a random number from a uniform distribution defined on the interval $[0,t_i)$ which is then recorded as the time of censoring.

%
\subsection{Retrieval accuracy}
%

It will be helpful to compare the retrieved $\matX^*$ with the `true' values $\matX$. For this purpose we choose the specific latent variables plotted in Figure \ref{special_X} (a) which are arranged in a specific geometrical pattern that allow quantitative measures of similarity to be defined. The samples that belong to either of the two circles, for instance, should be equidistant from the origin. If $\tilde{r}$ is the mean distance from the origin then we can define the radial error as
\begin{equation}
\mathcal{E}_{radial} = \frac{1}{|C|}\sum_{i\in C}\frac{|\vecx_i|-\tilde{r}}{\tilde{r}}
\end{equation}
where $C$ is the set of points belonging to the circle and $|C|$ is the number of samples belonging to that set. The error for both circles are averaged.

Similarly, the angles between each pair of samples belonging to each circle should be equal. In the case of the larger circle the angular separation should be $\tilde{\theta} = 2\pi/20$. If we let $\Delta\theta_i$ denote the angle between $\vecx_i$ and the neighbouring point then we can define the mean angular error as
\begin{equation}
\mathcal{E}_{angular} = \frac{1}{|C|}\sum_{i\in C}\frac{\Delta\theta_i-\tilde{\theta}}{\tilde{\theta}}.
\end{equation}

For both of the lines we can attempt a linear fit by writing $q_2 = \alpha q_1$. The value of $\alpha$ which minimises the sum of squared errors $\sum_i(x_{i2} - \alpha x_{i1})^2$ is given by $\hat{\alpha} = \sum x_{i1}x_{i2}/\sum x_{i1}^2$. We can then define the total sum of squares $SS_{tot} = \sum (x_{i2} - \left<x_{i2}\right>)^2$ and the sum of squared residuals $SS_{err} = \sum (x_{i2} - \hat{\alpha} x_{i1})^2$ and finally define
\begin{equation}
\mathcal{E}_{linear} = \frac{SS_{err}}{SS_{tot}}
\end{equation}
Note that $1 - \mathcal{E}_{linear}$ is called the coefficient of determination and is typically denoted by $R^2$ and takes a value between 0 and 1 where 1 corresponds to a perfect linear fit. These misalignment error measures have two desirable properties. Firstly, all three errors are zero for the `true' latent variables. Secondly, the error measures are invariant under rescaling of $\matX$.

%
\subsection{Retrieval accuracy of combined model compared to GPLVM}
%

\begin{figure}[h]
\centering
\includegraphics[scale = 0.8]{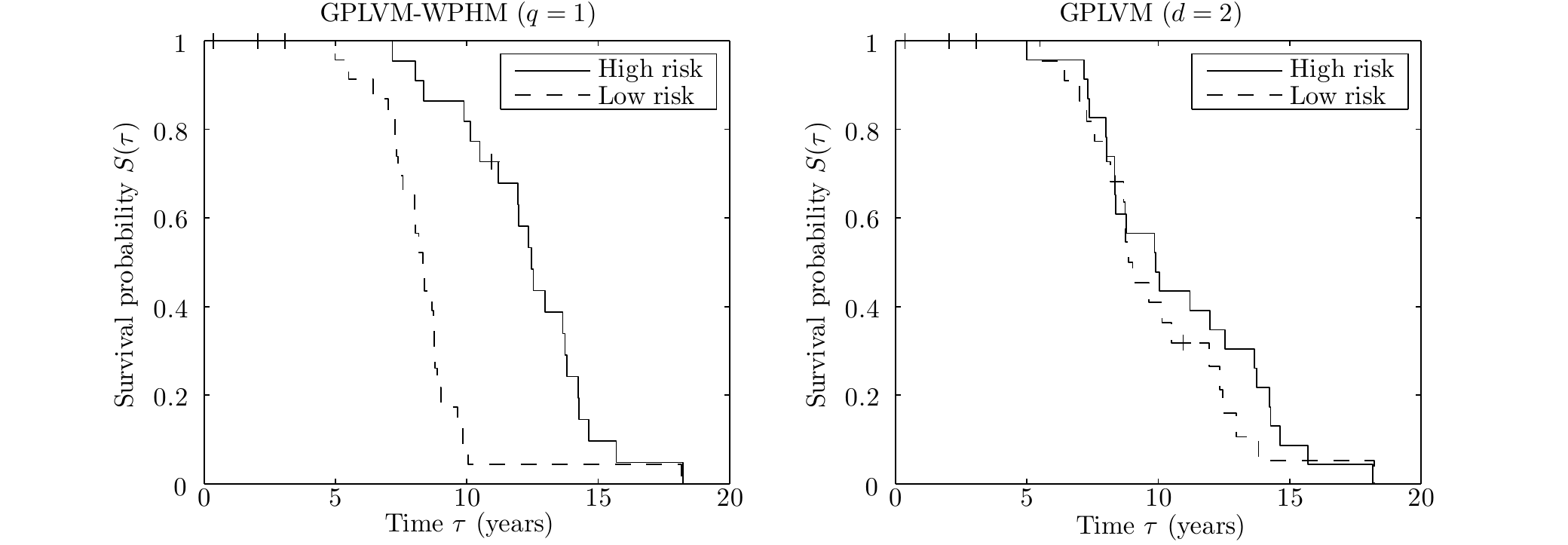}
\caption{Kaplan-Meier survival curves obtained in the latent variable space $q=1$ (left) and observed data space with $d=2$ (right) from synthetic two dimensional data that lie on a one dimensional non-linear manifold. A log rank test returned a p-value of 0.00006 on the left and 0.60755 on the right. Individuals were split into `high' and `low' risk groups by on the basis of risk factors $\vecb\cdot\vecx_i$. In the left figure we can clearly distinguish between the high and low risk individuals but this pattern is lost after the data are projected (in a non-linear manner) into the high dimensional space. A squared exponential kernel was used with $\beta^2=0.001$, $l=1$ and $\sigma=1$ to generate the two dimensional observed data. These results illustrate how the model can be used to extract a low dimensional manifold which may reveal additional structure to the data.}
\label{non-linear}
\end{figure}

Our aim here is to see if including survival data improves the ability of the model to accurately extract the correct low dimensional structure. To compare the combined GPLVM-WPHM to the GPLVM (which ignores the survival data) we generated fifty datasets from the two dimensional pattern in Figure \ref{special_X} (a). A linear kernel was used with $\beta=0.01$ and $d=10$. Survival times were generated as described above.

For each of the fifty datasets the GPLVM-WPHM was used to generate an optimal latent variable solution $\matX^*$ with $q=2$ and the misalignment errors were computed. The GPLVM was also used to generate a $q=2$ representation and misalignment errors were also computed for these solutions. Averaged over the fifty datasets a decrease was observed in the misalignment errors as shown in Table \ref{slvm:table:errordecrease}. We conclude that the inclusion of survival outcomes provides useful information that aids recovery of the true latent structure.

\begin{table}[b!]
  \begin{center}
    \begin{tabular}{|c|c|c|c|}
    \hline
    $\beta$ & $\mathcal{E}_{radial}$ & $\mathcal{E}_{angular}$ & $\mathcal{E}_{linear}$\\
	\hline
	0.1 &   $-7.3\%$   & $-5.3\%$  &  $-6.1\%$\\
    0.5 &    $-14.5\%$  &  $-17.1\%$  &  $-19.5\%$\\
    1.0 & $-16.0\%$ &  $-24.8\%$ & $-15.7\%$\\
    \hline
    \end{tabular}
  \end{center}
\caption[Percentage change in misalignment error when survival data are included]{The average percentage change in error when the GPLVM-WPHM is used instead of the GPLVM. A decrease in the misalignment error is observed due to the additional information provided by the survival data. The benefit becomes more apparent as the observed data become nosier (the survival data will contain roughly the same amount of `noise' in each experiment because they are generated with the same parameters throughout).}
\label{slvm:table:errordecrease}
\end{table}

%
\subsection{Integration of multiple sources}
%

Above we saw that including survival outcomes increases the accuracy of the retrieved latent variables. Now we investigate whether including multiple datasets simultaneously leads to similar improvement. We generated one dataset with $d_1=10$ and $\beta^2_1=0.1$ and second with $d_2=100$ and $\beta^2_2=1.0$ (a linear kernel was used). We compute the misalignment errors after analysing each dataset separately with the GPLVM-WPHM and compare this to the errors obtained after including both datasets simultaneously in the GPLVM-WPHM. The results in Table \ref{slvm:table:integration} show that it is beneficial to include both data sources together. These results were averaged over 50 repetitions.

\begin{table} 
  \begin{center}
    \begin{tabular}{|c|c|c|c|}
    \hline
     & $\mathcal{E}_{radial}$ & $\mathcal{E}_{angular}$ & $\mathcal{E}_{linear}$\\
	\hline
	$\matY_1$ ($d_1=10$, $\beta_1^2=0.1$) &   $0.0071$  &  $0.0093$  &  $0.0270$\\
    $\matY_2$ ($d_2=100$, $\beta_2^2=1.0$)&   $0.0244$ &   $0.0148$ &   $0.0509$\\
    $\matY_2\text{ \& }\matY_2$&    $0.0046$ &    $0.0052$  &  0.0146\\
    \hline
    \end{tabular}
  \end{center}
\caption[Misalignment error when two datasets are combined]{Misalignment errors (averaged over 50 repetitions) decrease when both datasets are combined simultaneously. It is beneficial to include all available information simultaneously rather than performing separate analyses on each dataset.}
\label{slvm:table:integration}
\end{table}

\begin{table}[h]
  \begin{center}
    \begin{tabular}{|c|c|c|c|c|}
    \hline
      $d=10$ & $d=25$ & $d=50$ & $d=100$\\
	\hline
	  $+1.2\%$ & $ +14.7\%$ & $+26.6\%$ & $+43.4\%$\\
    \hline
    \end{tabular}
  \end{center}
\caption{Percentage change in the mean squared error between predicted event times and reported event times in the high dimensional space compared the low dimensional space. The prediction error increases as the dimensional of the observed data grows. The noise was fixed at $\beta = 0.01$.}
\label{slvm:table:mse_d}
\end{table}

%
\subsection{Prediction accuracy using the latent variables}
%

In this section want to try and see the effect of overfitting due to high dimensionality. We generate datasets $\matY$ of different dimensions with $N=200$ individuals from a randomly generated matrix $\matX$ with $q=2$. Each dataset is split into a training and test set of equal size. In the high dimensional space we train a WPHM model on the training  individuals and then use the trained model to predict the event time for those individuals in the test set as described in Section \ref{predictions}. We then compute the mean square error (MSE) between the predicted and reported event times (censored individuals are excluded from the test set).

Next we run the GPLVM-WPHM on the same training data and use the trained model to firstly infer $\vecx^*$ from $\vecy^*$ for each individual in the test set and subsequently predict an event time. Again the MSE is computed and we can compare the MSE in the latent space to that obtained in the observed data space. In Table \ref{slvm:table:mse_d} we can see that the MSE increases as the dimension of the observed data increases. The results are averaged over fifty datasets.

Note that these data were generated with a linear kernel so the increase in MSE is not due to non-linearities induced during the generation of the synthetic data. Also the noise level is relatively low ($\beta^2=0.01$) so the observed data are only slightly corrupted with noise. We conclude that the increase in MSE is due to high dimension alone (rather than noise or non-linearities).

We also examined the effect that the noise level has (for fixed $d$). We can see from Table \ref{slvm:table:mse_beta} that the MSE in general increases with the noise level. The unusually large value for $\beta=0.5$ is due to an `outlier' (that is, one particularly bad prediction in the high dimensional space).

Finally, we investigated the behaviour at different levels of censoring. In Table \ref{slvm:table:mse_cens} we can see that higher levels of censoring lead to a degradation in predictive performance.

\begin{table}[t]
  \begin{center}
    \begin{tabular}{|c|c|c|c|c|}
    \hline
      $\beta=0.01$ & $\beta=0.1$ & $\beta=0.5$ & $\beta = 1.0$\\
	\hline
	  $+1.2\%$   & $+2.7\%$  &  $+38.0\%$ & $+5.67\%$ \\
    \hline
    \end{tabular}
  \end{center}
\caption{Percentage change in the mean squared error between predicted event times and the actual event times when computed in the high dimensional space compared to the low dimensional one for different noise values ($d=10$ in all cases). In general the error increases as the high dimensional data become more noisy. The large value at $\beta=0.5$ was due to one particularly poor prediction in the high dimensional space.}
\label{slvm:table:mse_beta}
\end{table}

%
\subsection{Non-linear dimensionality reduction}
%

In this section we investigate the effects that a non-linear mapping can induce. We used the squared exponential kernel to project latent variables with $q=1$ to $d=2$. The observed data would not be considered `high' dimensional but they now lie on a non-linear one dimensional manifold. In Figure \ref{non-linear} we compare survival curves in both spaces (obtained after training a GPLVM-WPHM and a WPHM respectively). The cohort was split into equally sized `high' and `low' risk groups by ranking all individuals according to the values of $\vecb\cdot\vecx_i$ and then separating them into two groups of the same size. A complete loss of structure is observed in the two dimensional space due to the non-linearities, whereas after we extract the one dimensional non-linear manifold we see a clear separation of the two groups. We can also compare the inferred hyperparameters to those that were used to generate the data. The generating hyperparameters in this case are $(\beta^2,\sigma,l,b,\rho,\nu)=(0.0010,1.00,1.00,-1.00,10.0,10.0)$ and the inferred values are $(\beta^2,\sigma,l,b,\rho,\nu)=(0.0006,1.23,1.11,-0.68,9.70,10.3)$.

This illustrates that the GPLVM-WPHM is useful not only for cases where $d>N$ but also cases where non-linear structures can be extracted from the covariates that may potentially reveal additional patterns of survival.

\begin{table}[b]
  \begin{center}
    \begin{tabular}{|c|c|c|c|c|}
    \hline
      $p=0.10$ & $p=0.25$ & $p=0.50$ & $p=0.75$\\
	\hline
	  $+36.3\%$   & $+60.9\%$  &  $+106.9\%$ & $+208.0\%$ \\
    \hline
    \end{tabular}
  \end{center}
\caption{Percentage change in the mean squared error between predicted event times and the actual event times when computed in the high dimensional space compared to the low dimensional one for different censoring levels. The variable $p$ is the fraction of individuals that were censored. Higher censoring leads to a loss in predictive accuracy. The dimension was $d=25$ with $\beta=1.0$.}
\label{slvm:table:mse_cens}
\end{table}

%
\subsection{Dimensionality detection}
%
\begin{figure}[t!]
\centering
\includegraphics[scale = 0.8]{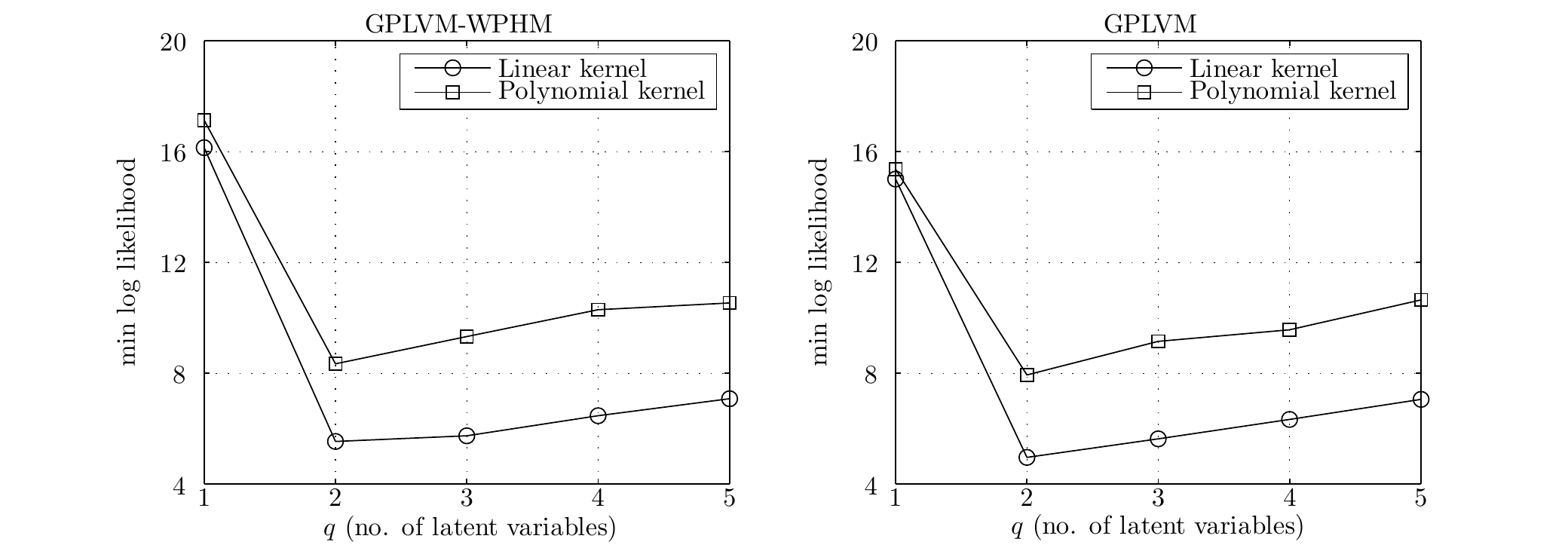}
\caption{Comparison of the GPLVM-WPHM and the GPLVM on synthetic data. The minimum negative log marginal likelihood is plotted as a function of $q$ for both models. Synthetic data with $d=10$, $\beta=0.1$ and $N=96$ (see Figure \ref{special_X} (a)) were generated using a linear kernel along with survival times. Both models correctly detect that $q^*=2$. However the model likelihood ratio between $q=2$ and $q=3$ using the GPLVM is 1.94 (that is, two latent variables is almost twice as probable as three), compared to a ratio of 1.23 using the GPLVM-WPHM. Furthermore, both models correctly detect that the linear kernel offers a better description of these data than the non-linear polynomial kernel.}
\label{dim}
\end{figure}

Next we illustrate the ability of the GPLVM-WPHM to correctly detect any intrinsic low dimensional structure. This is done by training models with different values of $q$ and comparing the minimum value of the negative log marginal likelihood. Shown in Figure \ref{dim} (a) is an example of the model correctly determining that $q^*=2$. Additionally we can compare this to an alternative kernel and we see that the linear kernel (correctly) offers the best description of these data.

In Figure \ref{dim} (b) we repeat the same experiment using the GPLVM and we see similar results. In fact the GPLVM has a slightly sharper minimum at $q=2$. One possible explanation for this is that the GPLVM-WPHM is overfitting slightly by using the third latent variable to explain some of the survival outcomes (the three regression coefficients are $b_1=-3.76$, $b_2 = 0.54$ and $b_3=1.19$).

%
\section{Experimental data}
%

Finally, we ran our model on a dataset of gene signature scores from breast cancer patients in the Guy's METABRIC dataset. This gene expression dataset from \cite{Curtis2012} was filtered for array intensity, quantile normalised, and batch-corrected for BeadChip (n=234 samples). In total there were $N=152$ patients and a total of $14,804$ gene expression levels per patient. We ranked the genes according to a univariate Cox model and selected the top 100 genes for use in our model. One of the aims of the METABRIC study was to search for potential gene signatures that are associated with clinical outcome (overall survival or progression free survival). Given the high dimensionality of the dataset there is a considerable risk of overfitting and therefore methods to probe associations between gene expression levels and survival outcomes while offering some protection against overfitting are needed. For this purpose we applied the GPLVM-WPHM to the dataset.

We performed 8-fold cross validation which consists of training a model on 7/8 of the data and then testing the predictive ability of the trained model on the remaining 1/8. As a metric of performance we used the C-statistic proposed by \cite{uno2011c} where larger values correspond to a greater ability to discriminate between high and low risk patients. The C-statistic was then averaged over the 8 cross validation runs.

In addition, we compared our model to three existing tools for high dimensional survival data. We performed L2 regularised Cox regression using the `penalised' R package \cite{GOE10}. Secondly, we used the `univariate shrinkage' method described in \cite{tibshirani2009univariate} and implemented in the `uniCox' R package. Finally, we used the `randomForestSRC' R package to implement a random forest for survival data \cite{r1,r2,r3}. 

\begin{figure}[t!]
\centering
\includegraphics[scale = 0.9]{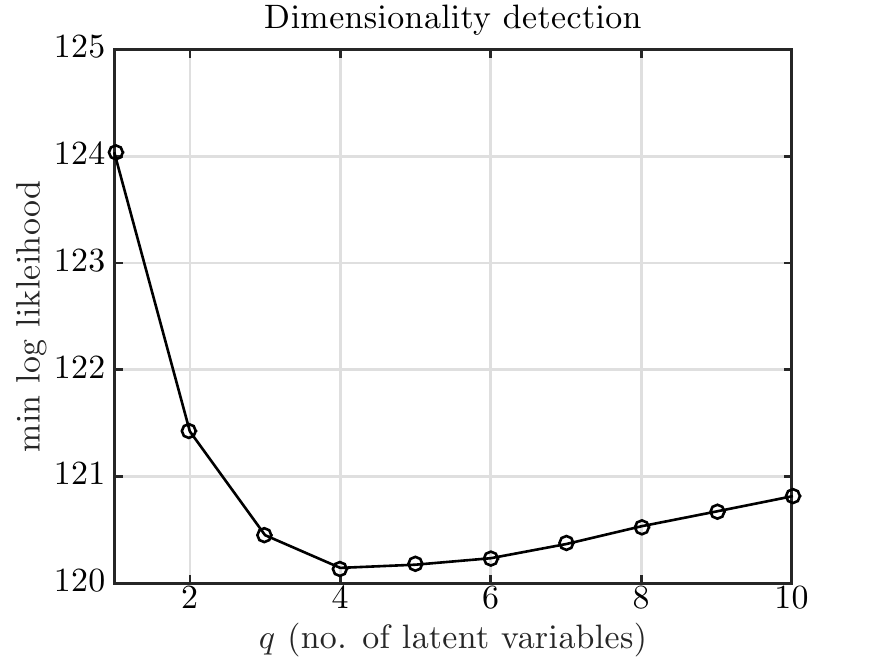}
\caption{Results from our analysis of experimental gene expression data. The minimum log likelihood corresponds to the most probably number of latent variables. In this case four latent variables are required to capture the information contained in 100 genes. Using this four dimensional representation leads to a considerable increase in predictive accuracy.}
\label{expt:metabric100}
\end{figure}

Using all $N=152$ of the samples we trained the GPLVM-WPHM for different values of $q$ and the minimum log likelihood obtained for models with $q=1,\ldots,10$ is plotted in Figure \ref{expt:metabric100}. The optimal number of latent variables is $q^*=4$ which indicates there is substantial redundancy between these genes.

The best performing models were the GPLVM-WPHM and the L2-regularised Cox model which both had an average C-statistic of 0.83. A WPHM model fitted to the $d=100$ genes had a score of 0.69. The univariate shrinkage method had a score of 0.76 while the random forest model performed the worst with a score of 0.70. There is therefore a significant improvement in predictive performance if a low dimensional representation is used instead of the original high dimensional data. The performance is superior or equal to the other models tested here. In addition, the information on intrinsic dimensionality may be of interest in itself, and the latent variable representation can be subsequently used for subsequent analyses.

\section{Discussion and conclusion}
\label{sec:disc}

Our proposed GPLVM-WPHM offers a novel supervised dimensionally reduction method for survival data. Simulation studies illustrate that including survival data is worthwhile and leads to more accurate retrieval of low dimensional structure. Our results also show that reducing the dimension can lead to a significant improvement in predictive accuracy as the effects of overfitting are diminished. In addition, our model can be used to extract non-linear low dimensional structure that has the potential to provide new insight into survival outcomes.

We used a gene expression dataset from the METABRIC study to show that using the GPLVM-WPHM can achieve a greater predictive accuracy than using the original dataset. This translates into a greater ability to discriminate between high and low risk breast cancer patients that may potentially be of practical benefit. Furthermore, the GPLVM-WPHM offers state of the art performance when compared to existing models for high dimensional survival data.

Future work could involve combining the GPLVM with more sophisticated survival analysis models. See \cite{LU08, MART11, VAN13} for examples of models that allow for flexible non-linear covariate effects in the hazard rate. Another research direction would be to apply some of the sparse GP regression techniques in order to reduce the computational burden.

\section{Acknowledgments}

This work was funded under the European Commission FP7 Imagint Project, EC Grant Agreement number 259881. We would like to thank Arnie Purushotham, Tony Ng, Anita Grigoriadis and in particular Katherine Lawler for their assistance in accessing and preparing the Guy's METABRIC gene expression data. We are grateful to the reviewers for their helpful suggestions and comments.

\bibliographystyle{unsrt}
\bibliography{refs}

\appendix

\section{Model Inference}
\label{app:sec:inf}

In Section 2.4 (Inference of parameters and hyperparameters) we wish to infer the values of the latent variables $\matX$ and WPHM parameters. We begin by defining the log likelihood as the negative log of [6] (square brackets denote equation numbers from the main text):
\begin{align}
\mathcal{L}(\matX,\vecb,\rho,\nu;\{\beta\},\tv) &= -\frac{1}{N}\log p(\matX,\vecb,\rho,\nu|\{\matY\},D,\tv,\{\beta\})\nonumber\\
&=\sum_{s=1}^S\left[\frac{d_s}{2N}\text{tr}(\matK_s^{-1}\matS_s)  +\frac{d_s}{2N}\log|\matK_s| + \frac{d_s}{2}\log 2\pi\right]-\frac{1}{N}\sum_{i:\Delta_i=1}\left[\log\lambda_0(\tau_i)+\vecb\cdot\vecx_i\right]\nonumber\\
&\qquad  +\frac{1}{N}\sum_{i=1}^N\Lambda_0(\tau_i)e^{\vecb\cdot\vecx_i} -\frac{1}{N}\log p(\matX)-\frac{1}{N}\log p(\vecb)-\frac{1}{N}\log p(\rho)-\frac{1}{N}\log p(\nu)\nonumber\\
&\qquad+const.
\label{eq:LL}
\end{align}
The marginal likelihood [7] is required for the hyperparameter posterior [10]. A Laplace approximation was constructed by expanding the log likelihood (\ref{eq:LL}) to second order about the minimum $(\matX^*,\vecb^*,\rho^*,\nu^*)$ which allows us to evaluate the integral in [7]. For compactness we write $\vecw = (\matX,\vecb,\rho,\nu)$:
\begin{align}
p(\{\matY\},D|\{\beta^2\},\tv) &= \int \text{d}\vecw\,e^{-N\mathcal{L}(\vecw)}\nonumber\\
&\approx\int\text{d}\vecw\, e^{-N\mathcal{L}(\vecw^*) - \frac{N}{2}(\vecw-\vecw^*)\cdot\matH(\vecw-\vecw^*)}\nonumber\\
&\propto p(\{\matY\},D|\matX^*,\vecb^*,\rho^*,\nu^*,\tv,\{\beta^2\})|N\matH^{-1}(\{\beta^2\},\tv)|^{1/2}.
\label{intro:eq:laplace}
\end{align}
The hessian matrix $\matH$ contains all of the second order partial derivatives
\begin{equation}
\matH(\{\beta^2\},\tv)=
\left(
\begin{array}{cc}
\matH_{XX} & \matH_{\tilde{b}X}  \\
\matH_{X\tilde{b}} & \matH_{\tilde{b}\tilde{b}} \\
\end{array}
\right)
\end{equation}
where for brevity we let $\tilde{\vecb} = (\vecb,\rho,\nu)$. The block matrices are defined by
\begin{equation}
[\matH_{XX}]_{p\eta,r\mu} = \frac{\partial^2\mathcal{L}}{\partial x_{p\eta}\partial x_{r\mu}}\text{,}\quad[\matH_{X\tilde{b}}]_{p\eta,\mu} = \frac{\partial^2\mathcal{L}}{\partial x_{p\eta}\partial \tilde{b}_{\mu}}\quad\text{and}\quad[\matH_{\tilde{b}\tilde{b}}]_{\eta,\mu} = \frac{\partial^2\mathcal{L}}{\partial \tilde{b}_{\eta}\partial\tilde{b}_{\mu}}.
\end{equation}
Second order partial derivatives are given below.

\subsection{Elimination of symmetries due to invariance under unitary transformations}

\begin{figure}[h]
\centering
\includegraphics[scale=0.75]{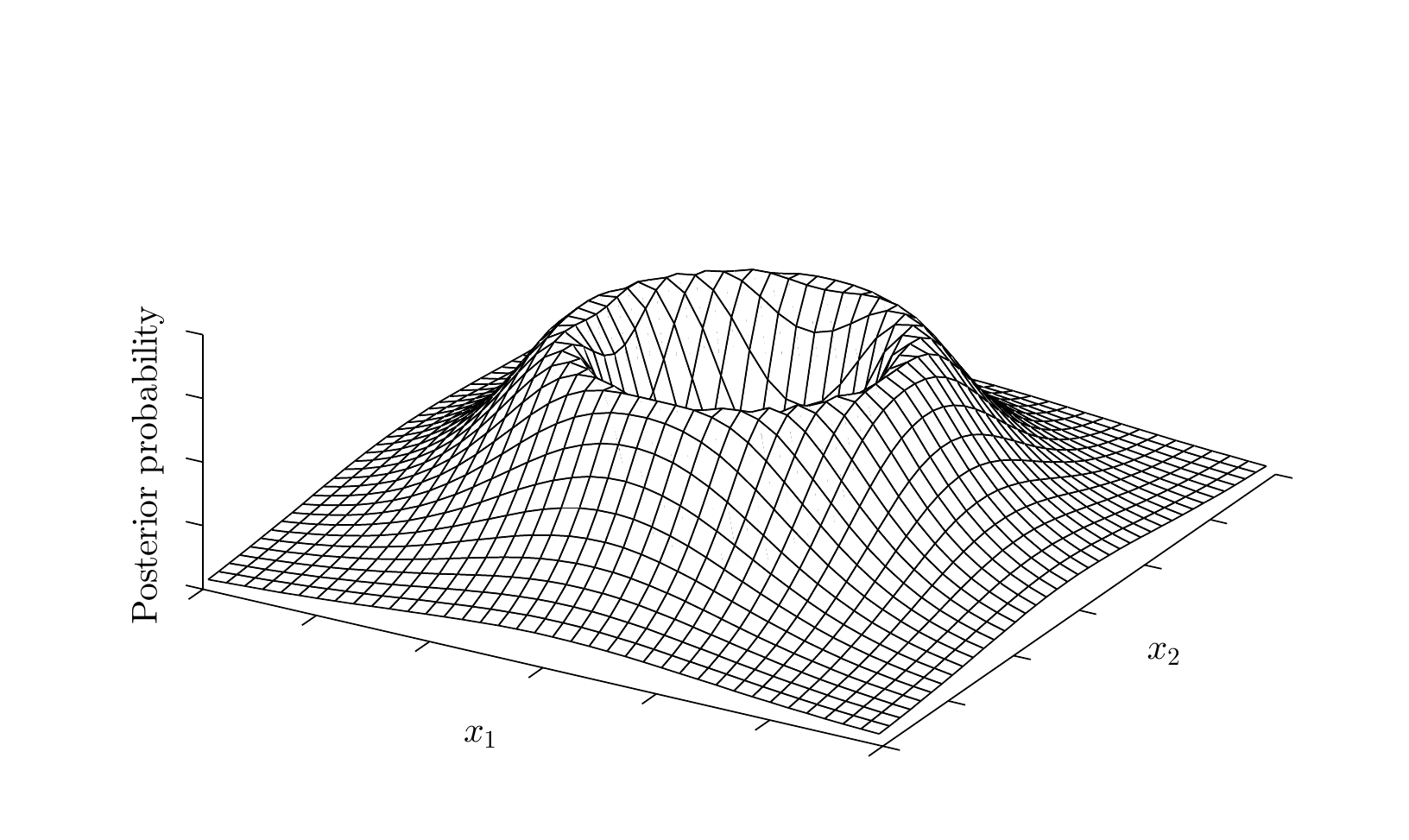}
\caption[Toy example of symmetries in the latent variable space]{Illustration of invariance under unitary transformations. The posterior probability density from a toy dataset with $N=1$ and $\vecx_1 = (x_1,x_2)$ is plotted (on an arbitrary scale). A vector of observed data $\vecy_1\in\mathbb{R}^5$ was randomly generated and the noise level was set to $\beta=0.01$. The fact that the data likelihood is invariant under rotations and reflections of $\vecx_1$ through the origin is readily apparent in this case. The points of maximum probability form a circle about the origin. The optimal value of $\vecx_1$ that is reported by the optimisation algorithm will depend on the initial conditions (which are generated randomly). If either $x_1$ or $x_2$ are close to zero then one of the second order derivatives will be close to zero. This renders a Gaussian approximation of the posterior invalid since the covariance matrix is no longer positive definite. Once the rotation and reflection symmetries have been eliminated the posterior reduces to a one dimensional unimodal distribution.}
\label{gplvm:fig:invariance}
\end{figure}

As mentioned in Section 2.5 a problem arises during the Laplace approximation due to fact that in the $Nq$-dimensional posterior search space of latent variables there exist directions in which the second order partial derivatives are zero. These directions point along lines where the log likelihood is constant. This occurs due to the invariance of the log likelihood to rotation or reflection of the latent variables. To see this let $\matU$ be a unitary matrix (corresponding to a rotation or reflection), such that $\matU^{\text{T}}\matU = \matU\matU^{\text{T}} = \matI$, and let $\tilde{\vecx} = \matU\vecx$. Then  
\begin{align}
\tilde{\vecx}_i\cdot\tilde{\vecx}_j &= \vecx_i\matU^{\text{T}}\matU\vecx_j = \vecx_i\cdot\vecx_j\text{ and}\label{gplvm:eq:inv1}\\
(\tilde{\vecx}_i-\tilde{\vecx}_j)^2 &= (\vecx_i - \vecx_j)\matU^{\text{T}}\matU(\vecx_i - \vecx_j) = (\vecx_i-\vecx_j)^2.
\label{gplvm:eq:inv2}
\end{align}
All of the kernels considered in this paper depend on the covariates solely through expressions of the form (\ref{gplvm:eq:inv1}, \ref{gplvm:eq:inv2}) and consequently are invariant under unitary transformations. Since the log likelihood depends on the latent variables via the kernel function, it too is invariant under unitary transformations. An example of this is shown in Figure \ref{gplvm:fig:invariance}.

There are two undesirable consequences of this property. Firstly, the second order partial derivatives of (\ref{eq:LL}) may evaluate to zero and hence $\matH$ will not be positive definite. Secondly, it means that there is not a unique latent variable representation of a dataset but rather an infinite mutually equivalent number corresponding to different rotations and reflections.

A computationally straightforward solution to this problem is to `pin down' the latent variable representation such that the symmetries are eliminated. Assuming that we are working in the standard basis $\{\vece_1,\ldots,\vece_q\}$ we demand that $\vecx_1$ is `pinned' to the $\vece_1$-axis which can always be achieved through an appropriate unitary transformation. We then require the second individual to be confined to the $\vece_1$--$\vece_2$ plane. This continues for the first $q-1$ individuals. In practice this implies that we simply populate the $(q-1)(q-2)/2$ elements in the upper right hand triangle of $\matX$ with zeros and optimise with respect to the remaining latent variables.
\begin{equation}
\matX = \left(
\begin{array}{cccc}
\tilde{x}_{11} & 0 & 0 & 0\\
\tilde{x}_{21} & \tilde{x}_{22} & 0 & 0 \\
\tilde{x}_{31} & \tilde{x}_{32} & \tilde{x}_{33} & 0 \\
\tilde{x}_{41} & \tilde{x}_{42} & \tilde{x}_{43} & \tilde{x}_{44} \\
\vdots & \vdots & \vdots & \vdots \\
\end{array}
\right).
\end{equation}

To eliminate reflection symmetries we require $\tilde{x}_{11} \geq 0, \tilde{x}_{22} \geq 0,\ldots,\tilde{x}_{qq}\geq0$. Reflection symmetries do not lead to a problem with zero second order partial derivatives but to obtain a unique solution it may be desirable to impose the above non-negativity conditions.

Note that the above solution may fail to guarantee a unique solution if $|\vecx_1| \approx 0$ since `pinning' a zero vector to the $\vece_1$-axis will not constrain the remaining latent variables. Furthermore, if $x_{22} \approx 0$ then reflection symmetry may not be broken. Although in those cases the solution may no longer be unique, the problem of zero partial derivatives will still be avoided.

\section{GPLVM partial derivatives}
\label{app:sec:gplvm}

In Section 2.4 (Inference of parameters and hyperparameters) we require the first and second order partial derivatives of the log likelihood for the purposes of gradient based optimisation and the construction of the Laplace approximation. Here we give partial derivatives for the GPLVM part of the log likelihood. The following identities are used \cite{COOKBOOK}:
\begin{align}
\frac{\partial |\matK|}{\partial \matK}&=|\matK|\matK^{-1}\label{app:eq:mc1}\\
\frac{\partial \text{tr}(\matA\matK^{-1}\matB)}{\partial \matK}&=-(\matK^{-1}\matB\matA\matK^{-1})^{\text{T}}\label{app:eq:mc2}\\
\frac{\matK_{kl}^{-1}}{\partial \matK_{ij}} &= \matK^{-1}_{ki}\matK^{-1}_{jl}\label{app:eq:mc3}
\end{align}
Define, for the purposes of this Section,
\begin{equation}
\mathcal{L}(\matX) =\sum_{s=1}^S\left[\frac{d_s}{2N}\text{tr}(\matK_s^{-1}\matS_s)  +\frac{d_s}{2N}\log|\matK_s| + \frac{d_s}{2}\log 2\pi\right].
\end{equation}
The following is true for any type of kernel function
\begin{equation}
\frac{\partial }{\partial \matX}\mathcal{L}(\matX)=\sum_{s=1}^S\sum_{i,j=1}^N\frac{	\partial \mathcal{L}}{\partial \matK_{ij}^s}\frac{\partial \matK_{ij}^s}{\partial \matX}
\label{app:gplvm:chainrule}
\end{equation}
where from (\ref{app:eq:mc1}, \ref{app:eq:mc2}) 
\begin{equation}
\frac{\partial \mathcal{L}}{\partial \matK^s}=-\frac{d_s}{2N}\matK_s^{-1}\matS_s\matK_s^{-1}+\frac{d_s}{2N}\matK_s^{-1}.
\label{app:gplvm:GLL}
\end{equation}
In what follows we drop the index $s$ for clarity and derive the partial derivatives for the linear, squared exponential and polynomial kernels.

\subsubsection*{The linear kernel}
The linear kernel is defined by (with $\sigma=1$)
\begin{equation}
\matK_{ij} = \vecx_i\cdot\vecx_j + \beta^2\delta_{ij}.
\end{equation}
First order partial derivatives are
\begin{equation*}
\frac{\partial}{\partial \matX}\mathcal{L}=-\frac{d}{N}\matK^{-1}\matS\matK^{-1}\matX+\frac{d}{N}\matK^{-1}\matX.
\end{equation*}
The second order partial derivatives are
\begin{align}
\frac{\partial}{\partial x_{p\nu}} (-\matK^{-1}\matS\matK^{-1}\matX)_{r\mu}=&-(\matK^{-1}\matS\matK^{-1})_{rp}\delta_{\mu\nu}+(\matK^{-1}\matS\matK^{-1}\matX)_{p\mu}(\matK^{-1}\matX)_{r\nu} +\nonumber\\
&+(\matK^{-1}\matS\matK^{-1}\matX)_{r\nu}(\matK^{-1}\matX)_{p\mu}+(\matK^{-1}\matS\matK^{-1})_{rp}(\matX^T\matK^{-1}\matX)_{\nu\mu} +\nonumber\\
&+(\matX^T\matK^{-1}\matS\matK^{-1}\matX)_{\nu\mu}(\matK^{-1})_{rp}
\end{align}
and
\begin{equation}
\frac{\partial}{\partial x_{p\nu}} (\matK^{-1})_{r\mu}= (\matK^{-1})_{rp}\delta_{\mu\nu}-(\matK^{-1}\matX)_{r\nu}(\matK^{-1}\matX)_{p\mu}-(\matX^T\matK^{-1}\matX)_{\nu\mu}(\matK^{-1})_{rp}.
\end{equation}

\subsubsection*{The squared exponential kernel}
The squared exponential kernel function is defined as
\begin{equation}
\matK_{ij} = \sigma e^{-\frac{l}{2}(\vecx_i-\vecx_j)^2} + \beta^2\delta_{ij}.
\end{equation}
Evaluation of (\ref{app:gplvm:chainrule}) takes $\mathcal{O}(N^2)$ operations to compute. However, for $i, j \text{ and } r$ distinct
\begin{align}
\frac{\partial \matK_{ij}}{\partial x_{r\mu}} &= 0\label{app:eq:se1}\\ 
\frac{\partial \matK_{ir}}{\partial x_{r\mu}} &= \frac{\partial \matK_{ri}}{\partial x_{r\mu}} = l \sigma(x_{i\mu}-x_{r\mu})e^{-\frac{l}{2}(\vecx_i - \vecx_r)^2} = l(x_{i\mu}-x_{r\mu})\matK_{ir}\label{app:eq:se2}
\end{align}
This allows us to compute (\ref{app:gplvm:chainrule}) with $\mathcal{O}(N)$ operations since from (\ref{app:gplvm:chainrule}) we can write
\begin{equation}
\frac{\partial \mathcal{L}}{\partial x_{r\mu}} = 2\sum_i\frac{\partial \mathcal{L}}{\partial \matK_{ir}}\frac{\partial \matK_{ir}}{\partial x_{r\mu}}.
\label{app:gplvm:firstse}
\end{equation}
Second order partial derivatives are obtained by differentiating (\ref{app:gplvm:firstse}) again to obtain
\begin{align}
\frac{\partial \mathcal{L}}{\partial x_{p\nu}\partial x_{r\mu}}&=2\sum_i^N\left\{\frac{\partial}{\partial x_{p\nu}}\left[-\frac{d}{2N}\matK^{-1}\matS\matK^{-1} + \frac{d}{2N}\matK^{-1}\right]_{ir}\left[l(x_{i\mu}-x_{r\mu})\matK_{ir}\right]\right.\nonumber\\
&\qquad+\left.\left[-\frac{d}{2N}\matK^{-1}\matS\matK^{-1} + \frac{d}{2N}\matK^{-1}\right]_{ir}\frac{\partial}{ \partial x_{p\nu}}\left[l(x_{i\mu}-x_{r\mu})\matK_{ir}\right]\right\}.
\label{app:gplvm:secondse}
\end{align}
On the first line there will be two terms. Beginning with the first term and using (\ref{app:eq:mc2}) and (\ref{app:eq:se1}, \ref{app:eq:se2}) and we can write
\begin{align}
\frac{\partial}{\partial x_{p\nu}}(\matK^{-1}\matS\matK^{-1})_{ir} &= \frac{\partial}{\partial x_{p\mu}}\sum_{tl}\matK^{-1}_{it}\matS_{tl}\matK^{-1}_{lr}\nonumber\\
& = \sum_{tl}\matK^{-1}_{it}\matS_{tl}\left[\frac{\partial}{\partial x_{p\nu}}\matK^{-1}_{lr}\right]+\sum_{tl}\left[\frac{\partial}{\partial x_{p\nu}}\matK^{-1}_{it}\right]\matS_{tl}\matK^{-1}_{lr}
\end{align}
which can be simplified to
\begin{align}
\frac{\partial}{\partial x_{p\mu}}(\matK^{-1}\matS\matK^{-1})_{ir} = &\sum_{k=1}^N\frac{\partial\matK_{pk}}{\partial x_{p\nu}}\left( -[\matK^{-1}\matS\matK^{-1}]_{ik}\matK^{-1}_{pr} -[\matK^{-1}\matS\matK^{-1}]_{ip}\matK^{-1}_{kr} \right.\nonumber\\
&\qquad\qquad\qquad\left.- [\matK^{-1}\matS\matK^{-1}]_{pr}\matK^{-1}_{ik}-[\matK^{-1}\matS\matK^{-1}]_{kr}\matK^{-1}_{ip}\right).
\label{app:gplvm:dinvKSinvK}
\end{align}
From (\ref{app:eq:mc3}) and (\ref{app:eq:se1}, \ref{app:eq:se2}) we can write the second term as
\begin{equation}
\frac{\partial}{\partial x_{p\nu}}\matK^{-1}_{ir} = \sum_{k=1}^N\frac{\partial \matK_{pk}}{\partial x_{p\nu}}\left(-\matK^{-1}_{ik}\matK^{-1}_{pr} -\matK^{-1}_{ip}\matK^{-1}_{kr}\right).
\label{app:gplvm:dinvK}
\end{equation}
Finally, the term on the second line of (\ref{app:gplvm:secondse}) is
\begin{align*}
\left[-\frac{d}{2N}\matK^{-1}\matS\matK^{-1} + \frac{d}{2N}\matK^{-1}\right]_{pr}\left[l^2(x_{p\mu}-x_{r\mu})(x_{r\nu}-x_{p\nu}) + l\delta_{\mu\nu}\right]\matK_{pr}&\qquad\text{when $i= p$ and $r \neq p$}\\
\left[-\frac{d}{2N}\matK^{-1}\matS\matK^{-1} + \frac{d}{2N}\matK^{-1}\right]_{ip}\left[l^2(x_{i\mu}-x_{p\mu})(x_{i\nu}-x_{p\nu}) - l\delta_{\mu\nu}\right]\matK_{ip}&\qquad\text{when $i\neq p$ and $r = p$}
\end{align*}
and 0 otherwise.

\subsubsection*{The polynomial kernel}
The polynomial kernel is defined by
\begin{align}
\matK_{ij} &= \sigma(1+\vecx_i\cdot\vecx_j)^{\alpha}+ \beta^2 \delta_{ij}\nonumber\\
&= \sigma\sum_{n=0}^{\alpha} \binom{\alpha}{n}(\vecx_i\cdot\vecx_j)^n + \beta^2 \delta_{ij}.
\end{align}
The partial derivatives of $\matK$ with respect to $\matX$ are
\begin{align}
\frac{\partial \matK_{ij}}{\partial x_{r\mu}} &= 0\\
\frac{\partial \matK_{ir}}{\partial x_{r\mu}}&= \frac{\partial \matK_{ri}}{\partial x_{r\mu}} = \sigma\sum_{n=1}^{\alpha} \binom{\alpha}{n}n(\vecx_i\cdot\vecx_r)^{n-1}x_{i\mu}\\
\frac{\partial \matK_{rr}}{\partial x_{r\mu}}&=2\sigma\sum_{n=1}^{\alpha} \binom{\alpha}{n}n(\vecx_r\cdot\vecx_r)^{n-1}x_{r\mu}.
\end{align}
Insertion into (\ref{app:gplvm:chainrule}) yields
\begin{equation}
\frac{\partial \mathcal{L}}{\partial x_{r\mu}} = 2\sum_i\frac{\partial \mathcal{L}}{\partial \matK_{ir}}\frac{\partial \matK_{ir}}{\partial x_{r\mu}}.
\label{app:gplvm:firstpoly}
\end{equation}
Second order partial derivatives
\begin{align}
\frac{\partial \mathcal{L}}{\partial x_{p\nu}\partial x_{r\mu}}& = 2\sum_i^N\left\{\frac{\partial}{\partial x_{p\nu}}\left[-\frac{d}{2N}\matK^{-1}\matS\matK^{-1} + \frac{d}{2N}\matK^{-1}\right]_{ir}\left[\sigma\sum_{n=1}^{\alpha}\binom{\alpha}{n}n(\vecx_i\cdot\vecx_r)^{n-1}x_{i\mu}\right]\right.\nonumber\\
&\qquad+\left.\left[-\frac{d}{2N}\matK^{-1}\matS\matK^{-1} + \frac{d}{2N}\matK^{-1}\right]_{ir}\frac{\partial}{ \partial x_{p\nu}}\left[\sigma\sum_{n=1}^{\alpha}\binom{\alpha}{n}n(\vecx_i\cdot\vecx_r)^{n-1}x_{i\mu}\right]\right\}.
\label{app:gplvm:secondpoly}
\end{align}
As in the case of the squared exponential kernel the first line contains two terms. They are 
\begin{align}
\frac{\partial}{\partial x_{p\nu}}(\matK^{-1}\matS\matK^{-1})_{ir} &= \frac{\partial}{\partial x_{p\mu}}\sum_{tl}\matK^{-1}_{it}\matS_{tl}\matK^{-1}_{lr}\nonumber\\
& = \sum_{tl}\matK^{-1}_{it}\matS_{tl}\left[\frac{\partial}{\partial x_{p\nu}}\matK^{-1}_{lr}\right]+\sum_{tl}\left[\frac{\partial}{\partial x_{p\nu}}\matK^{-1}_{it}\right]\matS_{tl}\matK^{-1}_{lr}
\end{align}
which can be simplified to (with $\alpha = 2$)
\begin{align}
\frac{\partial}{\partial x_{p\mu}}(\matK^{-1}\matS\matK^{-1})_{ir} = &\sum_{k=1}^N2\sigma x_{k\nu}(1+\vecx_k\vecx_p))\left( -[\matK^{-1}\matS\matK^{-1}]_{ik}\matK^{-1}_{pr} -[\matK^{-1}\matS\matK^{-1}]_{ip}\matK^{-1}_{kr} \right.\nonumber\\
&\qquad\qquad\qquad\left.- [\matK^{-1}\matS\matK^{-1}]_{pr}\matK^{-1}_{ik}-[\matK^{-1}\matS\matK^{-1}]_{kr}\matK^{-1}_{ip}\right).
\label{app:gplvm:dinvKSinvKpoly}
\end{align}
The second term is (with $\alpha = 2$)
\begin{equation}
\frac{\partial}{\partial x_{p\nu}}\matK^{-1}_{ir} = \sum_{k=1}^N\Bigg[-\matK^{-1}_{ik}\matK^{-1}_{pr} -\matK^{-1}_{ip}\matK^{-1}_{kr}\Bigg](2\sigma x_{k\nu}(1+\vecx_k\cdot\vecx_p)).
\label{app:gplvm:dinvKpoly}
\end{equation}
Terms from the second line of (\ref{app:gplvm:secondpoly}) are (with $\alpha = 2$)
\begin{align*}
\left[-\frac{d}{2N}\matK^{-1}\matS\matK^{-1} + \frac{d}{2N}\matK^{-1}\right]_{ip}2\sigma x_{i\mu}x_{i\nu}&\qquad\text{when $i\neq p$ and $r = p$}\\
\left[-\frac{d}{2N}\matK^{-1}\matS\matK^{-1} + \frac{d}{2N}\matK^{-1}\right]_{pp}2\sigma\left(\delta_{\mu\nu}(1+\vecx_p^2) + 2x_{p\mu}x_{p\nu}\right)&\qquad\text{when $i=p$ and $r = p$}\\
\left[-\frac{d}{2N}\matK^{-1}\matS\matK^{-1} + \frac{d}{2N}\matK^{-1}\right]_{pr}2\sigma\left(\delta_{\mu\nu}(1+\vecx_p\cdot\vecx_r) + x_{p\mu}x_{r\nu}\right)&\qquad\text{when $i=p$ and $r\neq p$}
\end{align*}
and zero otherwise.

\subsubsection*{Implementational details}

The sums over $k$ in (\ref{app:gplvm:dinvKSinvK}, \ref{app:gplvm:dinvK}) and (\ref{app:gplvm:dinvKSinvKpoly}, \ref{app:gplvm:dinvKpoly}) and the sums over $i$ in (\ref{app:gplvm:secondse}, \ref{app:gplvm:secondpoly}) can be eliminated by performing vectorised operations over appropriately defined matrices in \verb+matlab+. Matrices such as $\matX\matX^T$ and $\matK^{-1}\matS\matK^{-1}$ can be computed outside any loops. Since the Hessian matrix is symmetric it is necessary only to compute $Nq(Nq-1)/2$ partial derivatives.

\section{Weibull proportional hazards model (WPHM)}
\label{app:sec:wc}

In Section 2.4 (Inference of parameters and hyperparameters) we require partial derivatives of the log likelihood for the WPHM model. Define
\begin{align}
\mathcal{L}(\bv,\rho,\nu) &= -\frac{1}{N}\sum_{i:\Delta_i=1}\left[\log\lambda_0(\tau_i)+\bv\cdot\vecx_i\right] +\frac{1}{N}\sum_{i=1}^N\Lambda_0(\tau_i)e^{\bv\cdot\vecx_i}\nonumber\\
&\qquad\qquad -\frac{1}{N}\log p(\bv)-\frac{1}{N}\log p(\rho)-\frac{1}{N}\log p(\nu).
\end{align}
Partial derivatives are
\begin{equation}
\frac{\partial}{\partial \beta_s}\mathcal{L}(\bv,\rho,\nu) = -\frac{1}{N}\sum_{i:\Delta_i=1}x_{is} +\frac{1}{N}\sum_{i=1}^N \Lambda_0(\tau_i)x_{is}e^{\bv\cdot\vecx_i}
\label{app:wc:beta}
\end{equation}
and
\begin{align}
\frac{\partial}{\partial \rho}\mathcal{L}(\bv,\rho,\nu) &= \frac{N_1}{N}\frac{\nu}{\rho} + \frac{1}{N}\sum_{i=1}^N\frac{\partial\Lambda_0(\tau_i)}{\partial \rho}e^{\bv\cdot\vecx_i}\label{app:wc:rho}\\
\frac{\partial}{\partial \nu}\mathcal{L}(\bv,\rho,\nu) &= -\frac{N_1}{N}\frac{1}{\nu} -\frac{1}{N}\sum_{i:\Delta_i=1}\log(\tau_i/\rho) + \frac{1}{N}\sum_{i=1}^N\frac{\partial\Lambda_0(\tau_i)}{\partial \nu}e^{\bv\cdot\vecx_i}\label{app:wc:nu}
\end{align}
where we have used
\begin{align}
\frac{\partial}{\partial \rho} \log \lambda_0(\tau)&= -\frac{\nu}{\rho}\\
\frac{\partial}{\partial \nu} \log \lambda_0(\tau)&= \frac{1}{\nu}+\log(\tau/\rho)
\end{align}
and
\begin{align}
\frac{\partial\Lambda_0(\tau)}{\partial \rho} &= -\frac{\nu}{\rho}\left(\frac{\tau}{\rho}\right)^{\nu}\\
\frac{\partial\Lambda_0(\tau)}{\partial \nu} &= (\log\tau-\log\rho)\left(\frac{\tau}{\rho}\right)^{\nu}.
\end{align}
Since we require $\rho>0$ we write it in the form 
\begin{equation}
\rho = (1+\rho_{LB} + \exp(\tilde{\rho}))
\label{app:wc:par}
\end{equation}
where $\tilde{\rho}\in\mathbb{R}$ and $\rho_{LB}\geq0$ is a lower bound on $\rho$ that can be set manually. This formulation allows the use of unconstrained optimisation functions to be used. However the partial derivatives now become
\begin{equation}
\frac{\partial\mathcal{L}}{\partial\tilde{\rho}} = \frac{\partial\mathcal{L}}{\partial\rho}\frac{\partial\rho}{\partial\tilde{\rho}}\quad\text{with}\quad\frac{\partial\rho}{\partial\tilde{\rho}}=\frac{e^{\tilde{\rho}}}{1+e^{\tilde{\rho}}}.
\end{equation}
We also require $\nu>0$ and the same formulation is used.

\subsubsection*{Second order partial derivatives}

\begin{equation}
\frac{\partial^2}{\partial \beta_r\partial \beta_s}\mathcal{L}(\bv,\rho,\nu) =  \frac{1}{N}\sum_{i=1}^N\Lambda_0(\tau_i) x_{is}x_{ir}e^{\bv\cdot\vecx_i}
\end{equation}
and
\begin{align}
\frac{\partial^2}{\partial \rho^2}\mathcal{L}(\bv,\rho,\nu) &= -\frac{N_1}{N}\frac{\nu}{\rho^2} + \frac{1}{N}\sum_{i=1}^N\left[\frac{\nu(\nu+1)}{\rho^2}\left(\frac{\tau_i}{\rho}\right)^{\nu}\right]e^{\bv\cdot\vecx_i}\\
\frac{\partial^2}{\partial \nu^2}\mathcal{L}(\bv,\rho,\nu) &= \frac{N_1}{N}\frac{1}{\nu^2} + \frac{1}{N}\sum_{i=1}^N(\log\tau_i-\log\rho)^2\left(\frac{\tau_i}{\rho}\right)^{\nu}e^{\bv\cdot\vecx_i}.
\end{align}
Finally we require
\begin{align}
\frac{\partial^2}{\partial\nu\partial\rho}\mathcal{L}(\bv,\rho,\nu) &=\frac{\partial^2}{\partial\rho\partial\nu}\mathcal{L}(\bv,\rho,\nu) = \frac{N_1}{N}\frac{1}{\rho} -\frac{1}{N}\sum_{i=1}^N\left[\frac{\nu}{\rho}(\log\tau_i-\log\rho)\left(\frac{\tau_i}{\rho}\right)^{\nu}+\frac{1}{\rho}\left(\frac{\tau_i}{\rho}\right)^{\nu}\right]\\
\frac{\partial^2}{\partial\rho\partial\beta_s}\mathcal{L}(\bv,\rho,\nu) &= -\frac{1}{N}\frac{\nu}{\rho}\sum_{i=1}^N x_{is}\left(\frac{\tau_i}{\rho}\right)^{\nu}e^{\bv\cdot\vecx_i}\\
\frac{\partial^2}{\partial\nu\partial\beta_s}\mathcal{L}(\bv,\rho,\nu) &= \frac{1}{N}\sum_{i=1}^N (\log\tau_i - \log\rho)x_{is}\left(\frac{\tau_i}{\rho}\right)^{\nu}e^{\bv\cdot\vecx_i}
\end{align}
Since we have written the parameters $\rho$ and $\nu$ in the form (\ref{app:wc:par}) the second order partial derivatives are in practice given by
\begin{align}
\frac{\partial^2\mathcal{L}}{\partial\tilde{\rho}^2} &=\frac{\partial^2\mathcal{L}}{\partial\rho^2}\left(\frac{\partial\rho}{\partial\tilde{\rho}}\right)^2 + \frac{\partial\mathcal{L}}{\partial\rho}\frac{\partial^2\rho}{\partial\tilde{\rho}^2}\quad\text{with}\quad\frac{\partial^2\rho}{\partial\tilde{\rho}^2}=\frac{e^{\tilde{\rho}}}{(1+e^{\tilde{\rho}})^2},\\
\frac{\partial^2\mathcal{L}}{\partial\tilde{\rho}\partial\tilde{\nu}} &= \frac{\partial^2\mathcal{L}}{\partial\tilde{\nu}\partial\tilde{\rho}} = \frac{\partial^2\mathcal{L}}{\partial\rho\partial\nu}\frac{\partial\rho}{\partial\tilde{\rho}}\frac{\partial\nu}{\partial\tilde{\nu}},\\
\frac{\partial^2\mathcal{L}}{\partial\tilde{\rho}\partial\beta_s}&=\frac{\partial^2\mathcal{L}}{\partial\rho\partial\beta_s}\frac{\partial\rho}{\partial\tilde{\rho}}\text{ and }\frac{\partial^2\mathcal{L}}{\partial\tilde{\nu}\partial\beta_s} = \frac{\partial^2\mathcal{L}}{\partial\nu\partial\beta_s}\frac{\partial\nu}{\partial\tilde{\nu}}.
\end{align}

\subsubsection*{Prior terms}

When we assume the the prior distributions from Section 2.3 (The combined GPLVM-WPHM) we will need to include some additional terms in the first and second order partial derivatives. These are
\begin{align}
\frac{\partial}{\partial b_{\mu}}\left[-\frac{1}{N}\log p(\vecb)\right] & = \frac{1}{N\sigma_0^2}b_{\mu}\\
\frac{\partial}{\partial \nu}\left[-\frac{1}{N}\log p(\nu)\right] & = -\frac{\kappa_0-1}{N\nu}+ \frac{1}{N\alpha_0}\\
\frac{\partial}{\partial \rho}\left[-\frac{1}{N}\log p(\rho)\right] & = -\frac{\kappa_1-1}{N\rho}+ \frac{1}{N\alpha_1}
\end{align}
and
\begin{align}
\frac{\partial^2}{\partial b_{\mu} \partial b_{\nu}}\left[-\frac{1}{N}\log p(\vecb)\right] & = \delta_{\mu\nu}\frac{1}{N\sigma_0^2}\\
\frac{\partial^2}{\partial \nu^2}\left[-\frac{1}{N}\log p(\nu)\right] & = \frac{\kappa_0-1}{N\nu^2}\\
\frac{\partial^2}{\partial \rho^2}\left[-\frac{1}{N}\log p(\rho)\right] & = \frac{\kappa_1-1}{N\rho^2}
\end{align}

%
\section{GPLVM-WPHM}
\label{app:sec:slvm}
%

Define
\begin{equation}
\mathcal{L}(\vecb,\rho,\nu) = -\frac{1}{N}\sum_{i:\Delta_i=1}\left[\log\lambda_0(\tau_i)+\vecb\cdot\vecx_i\right] +\frac{1}{N}\sum_{i=1}^N\Lambda_0(\tau_i)e^{\vecb\cdot\vecx_i}.
\end{equation}
Then
\begin{equation}
\frac{\partial \mathcal{L}}{\partial x_{r\mu}} = \frac{1}{N}\left[-\delta_{1,\Delta_r}b_{\mu}+b_{\mu}\Lambda_0(\tau_r)e^{\vecb\cdot\vecx_r}\right]
\end{equation}
and
\begin{equation}
\frac{\partial^2 \mathcal{L}}{\partial x_{p\eta}\partial x_{r\mu}} = \frac{1}{N}\delta_{pr}b_{\mu}b_{\eta}\Lambda_0(\tau_r)e^{\vecb\cdot\vecx_r}.
\end{equation}
Using (\ref{app:wc:beta})
\begin{equation}
\frac{\partial \mathcal{L}}{\partial x_{r\mu}\partial b_{\eta}} = \frac{1}{N}\left\{\Lambda_0(\tau_r)x_{r\eta}b_{\mu}e^{\vecb\cdot\vecx_r} + \delta_{\mu\eta}\left[\Lambda_0(\tau_r)e^{\vecb\cdot\vecx_r}-\delta_{1,\Delta_r}\right]\right\}.
\end{equation}
Using (\ref{app:wc:rho}) and (\ref{app:wc:nu}) we obtain
\begin{align}
\frac{\partial \mathcal{L}}{\partial x_{r\mu}\partial \rho} & = -\frac{1}{N}\frac{\nu}{\rho}\left(\frac{\tau_r}{\rho}\right)^{\nu}b_{\mu}e^{\vecb\cdot\vecx_r}\\
\frac{\partial \mathcal{L}_D}{\partial x_{r\mu}\partial \nu} & = \frac{1}{N}(\log\tau_r-\log\rho)\left(\frac{\tau_r}{\rho}\right)^{\nu}b_{\mu}e^{\vecb\cdot\vecx_r}.
\end{align}
Due to the parameterisation (\ref{app:wc:par}) we will use the following partial derivatives in practice
\begin{equation}
\frac{\partial \mathcal{L}_D}{\partial x_{r\mu}\partial \tilde{\rho}} = \frac{\partial \mathcal{L}_D}{\partial x_{r\mu}\partial \rho}\frac{\partial\rho}{\partial\tilde{\rho}}\quad\text{with}\quad\frac{\partial\rho}{\partial\tilde{\rho}}=\frac{e^{\tilde{\rho}}}{1+e^{\tilde{\rho}}}.
\end{equation}

%
w\label{app:sec:pdgppred}
%

In Section 2.6 (Making predictions for new individuals) we require the partial derivatives of $p(\vecy^*|\vecx^*)$ in order to minimise the negative log likelihood using gradient based numerical methods. Firstly, we have
\begin{align}
\frac{\partial}{\partial x_{\mu}^*} \mathcal{L}(\vecx^*) &= -\frac{1}{N\kappa^3}\left(\frac{\partial\kappa}{\partial x_{\mu}^*}\right)\sum_{\nu=1}^d (y^*_{\nu}-m_{\nu})^2\nonumber\\
&\qquad\qquad-\frac{1}{N\kappa^2}\sum_{\nu=1}^2(y^*_{\nu}-m_{\nu})\left(\frac{\partial m_{\nu}}{\partial x_{\mu}^*}\right) + \frac{d}{N\kappa}\left(\frac{\partial\kappa}{\partial x_{\mu}^*}\right).
\end{align}
For the linear kernel
\begin{equation}
\frac{\partial m_{\nu}}{\partial x^*_{\mu}} = \vecx_{\mu}\cdot\matK^{-1}\vecy_{\nu}
\end{equation}
and
\begin{equation}
\frac{\partial\kappa}{\partial x^*_{\mu}}= 2x_{\mu}^* - \vecx_{\mu}\cdot\matK^{-1}\veck - \veck\cdot\matK^{-1}\vecx_{\mu}
\end{equation}
where $\vecx_{\mu}\in\mathbb{R}^{N\times 1}$ is the $\mu$th column of $\matX$. For the squared exponential kernel these are
\begin{equation}
\frac{\partial m_{\nu}}{\partial x^*_{\mu}} = -l\sum_{i=1}^N(x_{\mu}^*-x_{i\mu})k_i[\matK^{-1}\vecy_{\nu}]_i
\end{equation}
and
\begin{equation}
\frac{\partial\kappa}{\partial x^*_{\mu}}= 2l\sum_{i=1}^N(x_{\mu}^*-x_{i\mu})k_i[\matK^{-1}\veck]_i,
\end{equation}
where $k_i = [\veck(\vecx^*,\matX)]_i$. For the polynomial kernel these are
\begin{equation}
\frac{\partial m_{\nu}}{\partial x^*_{\mu}} = 2\sum_{i=1}^N(1+\vecx^*\cdot\vecx_i)x_{i\mu}[\matK^{-1}\vecy_{\nu}]_i
\end{equation}
and
\begin{equation}
\frac{\partial\kappa}{\partial x^*_{\mu}}= 4(1+\vecx^*\cdot\vecx^*)x_{\mu}^{\star} - 4\sum_{i=1}^N(1+\vecx^*\cdot\vecx_i)x_{i\mu}[\matK^{-1}\veck]_i.
\end{equation}

\end{document}